\documentclass[12pt]{amsart}
\usepackage{mathrsfs,amsmath,amsthm,amssymb}
\newtheorem{thm}{Theorem}[section]

\newtheorem{rem}[thm]{\bf{Remark}}

\numberwithin{equation}{section}
\newtheorem{prob}{Problem}

\begin{document}

\title[]{$p$-Groups for which each outer $p$-automorphism centralizes only $p$ elements}%
\author{Alireza Abdollahi}%
\address{Department of Mathematics, University of Isfahan, Isfahan 81746-73441, Iran \\ and School of Mathematics, Institute for Research in Fundamental Sciences (IPM), P.O.Box: 19395-5746, Tehran, Iran }%
\email{a.abdollahi@math.ui.ac.ir}%
\email{alireza\_abdollahi@yahoo.com}%
\author{S. Mohsen Ghoraishi}%
\address{School of Mathematics, Institute for Research in Fundamental Sciences (IPM), P.O.Box: 19395-5746, Tehran, Iran}%
\email{ghoraishi@gmail.com}%

\begin{abstract}
 An automorphism of a group is called outer if it is not an inner automorphism.
 Let $G$ be a finite $p$-group. Then for every outer $p$-automorphism $\phi$ of $G$ the
subgroup $C_G(\phi)=\{x\in G \;|\; x^\phi=x\}$ has order $p$ if and only if $G$ is of order  at most $p^2$.
\end{abstract}
\maketitle
\section{\bf Introduction}
 An automorphism of a group is called outer if it is not an inner automorphism. Let $p$ be any prime number. An automorphism of a group is called $p$-automorphism if its order is a power of $p$.  For any automorphism $\phi$ of a group $G$, $C_G(\phi)$ denotes the subgroup  $\{x\in G \;|\; x^\phi=x\}$.
Berkovich and  Janko  proposed the following problem in \cite[Problem 2008]{BZ}.
\begin{prob}\label{prob1} Study the $p$-groups $G$ such that for every outer $p$-automorphism $\phi$ of $G$ the
subgroup $C_G(\phi)$ has order $p$.
\end{prob}
Here we completely determine the structure of requested $p$-groups $G$ in Problem \ref{prob1}.
\begin{thm}\label{main}
  Let  $G$ be a finite $p$-group. For every outer $p$-automorphism $\phi$ of $G$ the
subgroup $C_G(\phi)$ has order $p$ if and only if $G$ is  of order at most  $p^2$.
\end{thm}

\section{\bf Preliminaries Results}
We use the following results in the proof of Theorem \ref{main}.
\begin{rem}\label{G-S} {\rm
By a famous  result of  Gasch\"utz \cite{Gas}, if $G$ is a finite $p$-group of order greater than $p$, then   $G$ admits
an outer $p$-automorphism.
Schmid \cite{S}  extended Gasch\" utz's result as follows: if $G$ is a finite
nonabelian $p$-group, then $G$ admits an outer $p$-automorphism $\phi$ such that the center $Z(G)$ of $G$ is contained in $C_G(\phi)$.
     The reader may pay attention to \cite{AG} to see  more recent results on the existence of noninner automorphism of order $p$ for finite non-abelian $p$-groups, a conjecture proposed by Y. Berkovich (see  Problem 4.13 of \cite{MK}).}
\end{rem}
For any  group $G$, we denote by $Aut^{\Phi}(G)$ the subgroup of all automorphisms of $G$ acting trivially  on the factor group $G/\Phi(G)$, where $\Phi(G)$ denotes the Frattini subgroup of $G$, the intersection of all maximal subgroups of $G$. By a well-known result of Burnside, $Aut^\Phi (G)$ is a $p$-group whenever $G$ is a finite $p$-group. Note that the inner automorphism group $Inn(G)$ of $G$ is contained in $Aut^\Phi(G)$.
\begin{rem}[{\cite[Theorem]{Mul}}]\label{Mul}
 {\rm Let $G$ be a finite $p$-group which is neither elementary abelian nor extraspecial. Then  $Aut^{\Phi}(G)$ properly contains $Inn(G)$.}
\end{rem}
Let $G$ be any group and $\phi$ is an automorphism of $G$. Let $N$ be a $\phi$-invariant subgroup of $G$; i.e., $N^\phi\subseteq N$. If $N$ is normal in $G$,  the map defined on $G/N$  by $xN\mapsto x^\phi N$ for all $x\in G$ is an automorphism of $G/N$. We will denote the latter map by $\overline{\phi}$.
\begin{rem}[{\cite[Lemma 2.12]{K}}]\label{C}
 {\rm   Suppose that  $\phi$ is an automorphism group of a finite group $G$ and $N$ is a
   normal $\phi$-invariant  subgroup. Then $|C_{G/N}(\overline{\phi})|\leq |C_G(\phi)|$.}
\end{rem}
\section{\bf Proof of  Theorem \ref{main}}
Assume that for every outer $p$-automorphism $\phi$ of $G$ the subgroup $C_G(\phi)$ has order $p$.

Let $V$ be an elementary abelian group of order $p^d$ and $d>2$. Suppose that $V=\langle v_1,\dots, v_d\rangle$. Then the map defined by $v_1\mapsto v_1v_2$ and $v_i\mapsto v_i$ for all $i>1$ can be extended to the automorphism $\phi$ of $V$ such that $|C_V(\phi)|=p^{d-1}>p$. The order of $\phi$ is $p$ and so it is an outer $p$-automorphism of $V$. Therefore, it follows that if $G$ is elementary abelian,   the order of $G$ is at most $p^2$.

Let $S$ be an extraspecial $p$-group of order $p^3$. Assume that  $p>2$. Suppose first that the exponent of $S$ is $p$. Then $S$ has a presentation as follows:
$$\langle x,y \;|\; x^p=y^p=1,[x,y]^y=[x,y]=[x,y]^x\rangle.$$
Now the map defined by $x\mapsto xy$ and $y\mapsto y$ determines the noninner automorphism $\alpha$
of order $p$ such that $\langle y, [x,y]\rangle \leq C_S(\alpha)$. To see the former claim, one may use  von Dyck's Theorem, as the $x^\alpha$ and $y^\alpha$ satisfy the same relations as $x$ and $y$ do, $\alpha$ can be extended to an endomorphism of $S$. Since $S=\langle x^\alpha,y^\alpha \rangle$, $\alpha$ is an epimorphism and since $S$ is finite, $\alpha$ is an automorphism of $S$.  \\
Now suppose that $S$ is of exponent $p^2$. Then $S$ has a presentation as follows:
$$\langle x,y \;|\; x^{p^2}=y^p=1,x^y=x^{1+p}\rangle.$$
The map defined by $x\mapsto xy$ and $y\mapsto y$ determines the noninner automorphism $\beta$
of order $p$ such that $\langle y, [x,y]\rangle \leq C_S(\beta)$. Showing that $\beta$ is an automorphism of $S$ is similar to that of $\alpha$, one may use the presentation of $S$ and observe that $x^\beta$ and $y^\beta$ satisfy corresponding relations as $x$ and $y$ do respectively.  \\
Now assume that $S=Q_8$ the quaternion group of order $8$ or $S=D_8$ the dihedral group of order $8$.
We know that $Q_8$ and $D_8$ have the following presentations:
$$Q_8=\langle x,y \;|\; x^4=1, x^2=y^2, x^y=x^{-1}\rangle, \;\;\; D_8=\langle x,y \;|\; x^4=y^2=1, x^y=x^{-1} \rangle.$$
The map defined on $S$ by $x\mapsto x$ and $y\mapsto xy$ can be extended to the noninner automorphism $\alpha$ of order  $4$ and
$\langle x\rangle\leq C_{S}(\alpha)$. \\ The following way to obtain such an automorphism $\alpha$ for $S\in \{D_8,Q_8\}$ is suggested by the referee.
Let $D$ be the semidihedral group of order $16$ and $C=\langle  c \rangle$
the cyclic subgroup of index $2$ in $D$. Both types of $S$ are subgroups of $D$ (see e.g., Theorem 1.2 of \cite{B}). Let $\bar{c}$ be the conjugation of $D$ by $c$. Then the fixed points of the
 restriction  of $\bar{c}$ to $S$ constitute the intersection $C \cap S$ of order $4$. Clearly, the restriction is a noninner automorphism of $S$.  \\

It follows that if $G$ is an extraspecial $p$-group, then $|G|>p^3$.
Thus $G$ is a central product of an extraspecial  group $A$ of order $p^3$ and another extraspecial group $B$. By previous paragraph, $A$ has an outer $p$-automorphism $\theta$ leaving  $Z(A)$  elementwise fixed. Now it is not hard to see that the map $\bar{\theta}$ on $G$ defined by $ab\mapsto a^\theta b$ for all $a\in A$ and $b\in B$ is an outer $p$-automorphism  fixing both $Z(A)$ and $B$ elementwise. This contradicts the assumption, since $|Z(A)B|>p$.

Now assume that $G$ is neither elementary abelian nor extraspecial.
By Remark \ref{Mul}, there exists some $\phi\in Aut^\Phi(G)\setminus Inn(G)$ so that  $|C_G(\phi)|=p$ by hypothesis.
It follows from Remark \ref{C} that  $|C_{G/\Phi(G)}(\overline{\phi})|\leq p$.
Thus $|G/\Phi(G)|=|C_{G/\Phi(G)}(\overline{\phi})|\leq p$.
This means that $G$ is a cyclic $p$-group. If $G=\langle a\rangle $ and $|a|=p^n>p^2$, then
 $\phi:\,a\mapsto a^{p^{n-1}+1}$ is an automorphism of order $p$. Now $\langle a^p\rangle \leq C_G(\phi)$, a contradiction. Thus $|G|=p^2$.
 The converse  obviously holds. This completes the proof. $\hfill \Box$\\

\begin{center}{\textbf{Acknowledgments}}
\end{center}
The authors are grateful to the referee for his careful readings and pointing some inaccuracies in the first draft of the current paper.
The first author was financially supported by the Center of Excellence for Mathematics, University of Isfahan.
This research was in part supported by a grant IPM (No. 91050219).

\end{document}